\documentclass[reqno]{amsart}
\usepackage{hyperref}

\usepackage[centertags]{amsmath}
\usepackage{amsfonts}
\usepackage{amssymb}
\usepackage{amsthm}
\usepackage{newlfont}
\usepackage{fancyhdr}

\makeatletter
\@namedef{subjclassname@2020}{%
  \textup{2020} Mathematics Subject Classification}
\makeatother

\fancypagestyle{plain}{\fancyfoot[C]{\vspace*{2\baselineskip}\thepage}}%Change number to move page number
\theoremstyle{plain}
\newtheorem{theo}{Theorem}[section]

\newtheorem{lem}{Lemma}[section]

\newtheorem*{conj*}{Denjoy's Conjecture}

\begin{document}
\title[%EJDE-2012/177
\hfil Subnormal solutions of delay differential equations]
{Subnormal transcendental meromorphic solutions of difference equations with Schwarzian derivative}

\author[M. T. Xia, J. R. Long, X. X. Xiang \hfil %EJDE-2012/177\hfilneg
]
{ Mengting Xia, Jianren Long$^{*}$, Xuxu Xiang}

\address{Mengting Xia \newline
School of Mathematical Sciences, Guizhou Normal University, Guiyang, 550025, P.R. China.}
\email{2190331341@qq.com}

\address{Jianren Long \newline
School of Mathematical Sciences, Guizhou  Normal University, Guiyang, 550025, P.R. China.}
\email{longjianren2004@163.com }

\address{Xuxu Xiang \newline
School of Mathematical Sciences, Guizhou  Normal University, Guiyang, 550025, P.R. China.}
\email{1245410002@qq.com }

%\address{Ke-e Qiu \newline
%School of Mathematics and Computer Science, Guizhou Normal
%Colleage, Guiyang, 550018, P.R. China.}
%\email{qke456@sina.com }

\thanks{The research  was supported by the National Natural Science Foundation of China (Grant No. 12261023, 11861023).\\$*$Corresponding author.}
\subjclass[2020]{Primary 30D35; Secondary 34M10, 34M55}
\keywords{Delay differential equations; Subnormal solutions; Existence; Schwarzian derivative;  Painlev\'{e} equations.}

\begin{abstract}
The existence of subnormal solutions of following three difference equations with   Schwarzian  derivative 
$$\omega(z+1)-\omega(z-1)+a(z)(S(\omega,z))^n=R(z,\omega(z)),$$
$$\omega(z+1)\omega(z-1)+a(z)S(\omega,z)=R(z,\omega(z)),$$
and
$$(\omega(z)\omega(z+1)-1)(\omega(z)\omega(z-1)-1)+a(z)S(\omega,z)=R(z,\omega(z))$$
are studied by using Nevanlinna theory, where $n\ge 1$ is an integer, $a(z)$ is small with respect to $\omega$, $S(\omega,z)$ is Schwarzian derivative, $R(z,\omega)$ is rational in $\omega$ with small meromorphic coefficients with respect to $\omega$. The necessary conditions for the existence of subnormal transcendental meromorphic solutions of the above equations are obtained. Some examples are given to support these results.

\end{abstract}

\maketitle
\numberwithin{equation}{section}
\newtheorem{theorem}{Theorem}[section]
\newtheorem{lemma}[theorem]{Lemma}
\newtheorem{definition}[theorem]{Definition}
\newtheorem{example}[theo]{Example}
\newtheorem{remark}[theo]{Remark}
\allowdisplaybreaks

\section{Introduction}

In what follows, we assume the reader is familiar with the basic notations  of Nevanlinna theory, such as the characteristic function $T(r,\omega)$, proximity function $m(r,\omega)$, and counting function $N(r,\omega)$, where $\omega$ is a meromorphic function, see \cite{c, f} for more details. Let $\rho_2(\omega)$ denotes the hyper order of $\omega$. As usual, we use $S(r,\omega)$ to denote any quantity satisfying $S(r,\omega)= o(T(r, \omega))$ as $r$ tends to infinity, possibly outside an exceptional set of finite logarithmic measure. For a meromorphic function $g$, if $T(r,g)=S(r,\omega)$, we say $g$ is small with respect to $\omega$. Furthermore, a transcendental meromorphic function $\omega$ is called subnormal if it satisfies
$\displaystyle\limsup_{r\rightarrow\infty}\frac{\log T(r,\omega)}{r}=0$.

The classical Malmquist theorem \cite{M} implies that if the first order  differential equation
\begin{align}  
    \label{01}
\omega'(z)=R(z,\omega)
\end{align}
admits a transcendental meromorphic solution, where $R(z,\omega)$ is a rational function in $z$ and $\omega$, then \eqref{01} reduces to a differential Riccati equation.  Yosida \cite{Yosida}  and Laine \cite{Laine 1}  given elegant alternate proofs of the classical
 Malmquist theorem by using  Nevanlinna theory. 
 A precise classification of the differential equation 
 \begin{equation}
    \label{0}
    (\omega'(z))^n=R(z,\omega)
 \end{equation}
is given by  Steinmetz \cite{Steinmetz}, and Bank and Kaufman \cite{Bank}, where $n$ is a positive integer and $R(z,\omega)$ is rational in both arguments.   
See also \cite[Chapter 10]{Laine} for Malmquist–Yosida–Steinmetz type theorems. 

The Schwarzian derivative of a meromorphic function $\omega$ is defined as 
   $ S(\omega,z)=(\frac{\omega''}{\omega'})'-\frac{1}{2}(\frac{\omega''}{\omega'})^2=\frac{\omega'''}{\omega'}-\frac{3}{2}(\frac{\omega''}{\omega'})^2.$
The Schwarzian derivative plays a significant role in multiple branches of complex analysis \cite{EH,OL,NS}, particularly in the theories of univalent functions and conformal mappings. Research has further demonstrated profound connections between this operator and both second order linear differential equations \cite{Laine} and the Lax pairs of certain integrable partial differential equations \cite{JW}. In particular, the equation \eqref{0} can be rewritten  as $(\frac{\omega'(z)}{\omega(z)})^n=\frac{R(z,\omega)}{\omega^n(z)}=R_1(z,\omega)$, then by replace the logarithmic derivative $\frac{\omega'(z)}{\omega(z)}$ with the Schwarzian derivative $S(\omega,z)$,  Ishizaki \cite{I1991} established several Malmquist-type theorems for the equation the Schwarzian differential equation
\begin{align}
    \label{00}
    (S(\omega,z))^n=R(z,\omega),
\end{align}
where $n$ is a positive integer, and $R(z,\omega)$ is an irreducible rational function in $\omega$ with meromorphic coefficients. For equation \eqref{00} with polynomial coefficients, Liao and Ye \cite{Liao1} investigated the growth of meromorphic solutions. Recently, all transcendental meromorphic solutions of the autonomous  Schwarzian differential equations have been constructed in \cite{Liao,Zhang}.

 %Yanagihara \cite[Corollary 6]{Y} proved that  if the first order difference equation $\omega(z+1)=R(z,\omega)$ has nontrivial meromorphic solutions $\omega$ with $\rho_2(\omega)<1$, then this equation   must reduce into the difference Riccati equation, where  $R(z,\omega)$ is rational in both arguments. This is a natural difference analogue of Malmquist’s theorem on first-order differential equation \eqref{01}.  A natural difference
% analogue of \eqref{0}, i.e., the first-order difference equation $(\omega(z+1))^n=R(z,\omega)$ was studied by Zhang et al.\cite{zhangyy}.  

The second order differential equation   $\omega''=R(z,\omega)$
have been classified by Fuchs \cite{Fu},  Gambier \cite{Gambier} and Painlev\'{e} \cite{p1,p2},  and they obtained six equations, known as the  Painlev\'{e}  equations.
 Similar the fact that the second order differential equation can be reduced into  Painlev\'{e}  equations, Halburd et al. \cite{HK1,HK2} proved if   the second order difference equation
 \begin{align}
    \label{03}
    \omega(z+1)*\omega(z-1)=R(z,\omega)
 \end{align}
existences finite order meromorphic solutions, then this difference equation reducing into a short list of canonical equations, including the difference Painlev\'{e} I-III equations, where  operation $*$ stands either for the
 addition or the multiplication, $R(z, \omega)$ is rational in $\omega$ with small functions of $\omega$ as coefficients. Later, Ronkainen \cite{R}  singled out a class of equations containing  the difference Painlev\'{e} \uppercase\expandafter{\romannumeral5} equation
 from the  Painlev\'{e} \uppercase\expandafter{\romannumeral5} type difference equation 
\begin{align}
    \label{009}
    (\omega(z)\omega(z+1)-1)(\omega(z)\omega(z-1)-1)=R(z,\omega). 
\end{align}
The discrete (or difference) Painlev\'{e} equations were attracted by different researchers, for example, see \cite{AHH,HK1,w} and therein references.

It is also worth noting that reductions of integrable differential-difference equations may give rise to
 delay differential equations with formal continuum limits to  Painlev\'e equations. In \cite{QCS}, Quispel, Capel and Sahadevan shown the equation 
\begin{equation}
    \label{9.1}
    \omega(z)[\omega(z+1)-\omega(z-1)]+a\omega'(z)=b\omega(z),
\end{equation}
where $a$ and $b$ are constants, can be obtained from the symmetry reduction of the Kac-Van Moerbeke equation and has a formal continuous limit to the first Painlev\'{e} equation
$y''=6y^2+t.$ In 2017, Halburd and Korhonen \cite{HK} considered an extended version of \eqref{9.1} and obtained the following results.
\begin{theo} \cite{HK} \label{thmx1.1}
    Let $\omega$ be a transcendental meromorphic solution of  
    \begin{equation}
    \label{9.2}
    \omega(z+1)-\omega(z-1)+a(z)\frac{\omega'(z)}{\omega(z)}=R(z,\omega(z))=\frac{P(z,\omega(z))}{Q(z,\omega(z))},
    \end{equation}
    where $a(z)$ is rational in $z$, $P(z,\omega)$ is a polynomial in $\omega$ having rational coefficients in $z$, and $Q(z,\omega)$ is a polynomial in $\omega$ with $\deg_{\omega}(Q)\ge 1$ and the roots of $Q(z,\omega)$ are non-zero raional functions of $z$ and not the roots of $P(z,\omega)$. If $\rho_2(\omega)<1$, then
    $$\deg_\omega(Q)+ 1=\deg_\omega(P)\le 3$$ 
    or $\deg_{\omega}R(z,\omega)\le 1$.
\end{theo}

Recently, related results on the Theorem \ref{thmx1.1} have been obtained in \cite{cc,cc1,ls,lx,x2025,Zhang}.
Liu et al. \cite{ls} considered the delay differential equation  \eqref{9.2} also can be viewed as a combination of  second order difference equation  with  Malmquist-Yosida type differential equation \eqref{01}. Therefore,
by  combining the equation \eqref{03} with Malmquist-Yosida type differential equation \eqref{01}, Liu et al.  \cite{ls} considered delay differential 
\begin{equation}
    \label{1.2}
    \omega(z+1)\omega(z-1)+a(z)\frac{\omega'(z)}{\omega(z)}=R(z,\omega(z))=\frac{P(z,\omega(z))}{Q(z,\omega(z))},
    \end{equation}
and obtained the same results of Theorem \ref{thmx1.1}. Inspired by this idea, Du et al. \cite{dz} also considered the necessary conditions for the existence of
transcendental meromorphic solutions of
\begin{align}
    \label{012}
    (\omega(z)\omega(z+1)-1)(\omega(z)\omega(z-1)-1)+a(z)\frac{\omega'(z)}{\omega(z)}=R(z,\omega(z))=\frac{P(z,\omega(z))}{Q(z,\omega(z))},
\end{align}
 which is  a combination of Painlev\'{e} \uppercase\expandafter{\romannumeral5} type  difference equation  \eqref{009} with Malmquist-Yosida type differential equation \eqref{01}.

According to the discussion above and  inspired by the above works of Ishizaki \cite{I1991} and Liu et al.  \cite{ls},   the following questions naturally arise: 
 What would occur when replacing the logarithmic derivative term $\frac{\omega'(z)}{\omega(z)}$ in \eqref{9.2}, \eqref{1.2} and \eqref{012} with the Schwarzian derivative $S(\omega,z)$ respectively? What would occur when combining the difference equation \eqref{03} or \eqref{009} with Schwarzian differential equation \eqref{00}?  These questions prompt us to consider the following three difference equations:
\begin{equation}
    \label{9.5}
   \omega(z+1)-\omega(z-1)+a(z)(S(\omega,z))^n=R(z,\omega(z))=\frac{P(z,\omega(z))}{Q(z,\omega(z))},
\end{equation}
\begin{equation}
        \label{eq.1.1}
        \omega(z+1)\omega(z-1)+a(z)S(\omega,z)=R(z,\omega(z))=\frac{P(z,\omega(z))}{Q(z,\omega(z))},
    \end{equation}
    and 
\begin{equation}
    \label{eq.1.3}
        (\omega(z)\omega(z+1)-1)(\omega(z)\omega(z-1)-1)+a(z)S(\omega,z)=R(z,\omega(z))=\frac{P(z,\omega(z))}{Q(z,\omega(z))},
    \end{equation}
  where $n\ge 1$ is an integer, $a(z)$ is small with respect to $\omega$, $R(z,\omega)$ is rational in $\omega$ with small meromorphic coefficients with respect to $\omega$, $P(z,\omega)$, $Q(z,\omega)$ are polynomials in $\omega$ having meromorphic coefficients small with respect to $\omega$ in $z$.

In 1978, Bishop \cite{B} investigated soliton behavior in discrete nonlinear lattices, laying the foundation for numerical analysis of the Discrete Sine-Gordon equation (DSG). In 2019, Khare et al.\cite{KS} studied discrete Sine-Gordon solutions under nonuniform coupling. Later, Kevrekidis \cite{K}  employed neural networks to solve DSG soliton dynamics. In 2022, Pelinovsky \cite{PS} analyzed the stability of two-dimensional DSG systems, providing theoretical support for experiments involving optical vortex solitons.

When taking specific forms of $R(z,\omega)$ and $a(z)$, the equation \eqref{eq.1.1} can be reduced to a special case of the DSG equation. By setting $\omega(z)=\log(u_n)$ and appropriately choosing $R(z,\omega)$ and $a(z)$, the original equation can be transformed into the exponential form of the DSG equation. In this case:

\begin{itemize}
\item The product term $\omega(z+1)\omega(z-1)$ corresponds to the cross-nonlinear term in the DSG equation.
\item The Schwarzian derivative $S(\omega,z)$ reflects discrete curvature correction, similar to the derivative term in the continuous Sine-Gordon equation.
\end{itemize}

In 1981, Hirota \cite{HR} introduced the bilinear form of the discrete Korteweg-de Vries (KdV) equation, laying the foundation for constructing soliton solutions and auto-Bäcklund transformations. Later, Nijhoff and Capel \cite{NC} studied the discrete KdV equation as a reduction of the discrete KP hierarchy and derived auto-Bäcklund transformations using lattice integrability. In 2015, Tongas and Tsoubelis \cite{TT} investigated noncommutative generalizations and quantum deformations of the discrete KdV equation, thereby obtaining new auto-Bäcklund transformations.

The bilinear form of the equation \eqref{eq.1.3} is naturally connected to the Darboux transformation of discrete integrable systems. When $R(z,\omega)$ takes specific rational functions, the original equation is equivalent to the auto-Bäcklund transformation of the discrete KdV equation:

\begin{itemize}
\item The bilinear term $(\omega\omega_{\pm 1}-1)$ corresponds to the three-point relation of discrete KdV equation.
\item The Schwarzian derivative term $S(\omega,z)$ provides spectral parameter dependence.
\end{itemize}

These correspondences demonstrate that the difference equations with Schwarzian derivatives studied in this paper provide a universal framework for discrete integrable systems, fostering cross-fertilization among various fields in mathematical physics.
  
This paper is organized as follows. The existence of subnormal transcendental meromorphic solutions of  \eqref{9.5}, \eqref{eq.1.1} and \eqref{eq.1.3} are characterized in Sections 2 and 3, respectively. Some auxiliary results which can be used in the proof of our results are shown in Section 4. Main results are proved in Sections 5 and 6 respectively.

%Firstly, we assume that $\omega$ has finitely many poles. By \eqref{2.1}, we can get $N(r,S(\omega,z))=S(r,\omega)$. Then by Lemma \ref{lem1}, Lemma \ref{lem2} and Lemma \ref{lem3}, we obtain
%\begin{eqnarray*}
   % \deg_{\omega }(P)T(r,\omega )+S(r,\omega)&=&T(r,P(z,\omega))    \\
   % &=&T(r,\omega(z+1)\omega(z-1)+aS(\omega,z))    \\
   % &\le& 2T(r,\omega)+m(r,S(\omega,z))+N(r,S(\omega,z))+S(r,\omega)   \\
   % &\le& 2T(r,\omega)+S(r,\omega),
%\end{eqnarray*}
%which implies that $\deg_{\omega }(P)\le 2$.

\section{Second order difference equation with Schwarzian derivative}
%Recently, related results on the Theorem \ref{thmx1.1} have been obtained in\cite{cc,cc1,ls,lx,x2025,Zhang}.  In particular,motivated by the work of Ishizaki\cite{I1991},  Nie et al.\cite{n2025} considered a variant of \eqref{9.2} by replacing the logarithmic derivative $\frac{\omega'(z)}{\omega(z)}$  with the Schwarzian derivative $S(\omega,z)$. \begin{theo} \cite{n2025} \label{thmx1.2} Let $\omega$ be a subnormal transcendental meromorphic solution of the equation \begin{equation} \label{9.3}\omega(z+1)-\omega(z-1)+a(z)S(\omega,z)=\frac{P(z,\omega(z))}{Q(z,\omega(z))},  \end{equation}where $a(z)$ is rational, $P(z,\omega)$ and $Q(z,\omega)$ are coprime polynomials in $\omega$ with rational coefficients, then $\deg_\omega(R)\le 7$ and $\deg_\omega(P)\le \deg_\omega(Q)+1$. Moreover, if $Q(z,\omega)$ has a rational function root $b_1$ in $\omega$ with multiplicity $k$, then $k\le 2$.\end{theo}Also inspired by the works of Ishizaki\cite{I1991}, Theorems \ref{thmx1.2} and \ref{thmx1.3}, we considered the equation \eqref{9.5}, and get the following result.
Theorem \ref{the9} below shows the  necessary conditions for the existence of
subnormal transcendental meromorphic solutions of the \eqref{9.5}, which is a generalization of \cite[Theorem 1.1]{n2025}.
\begin{theo}  \label{the9} 
    Let $\omega$ be a subnormal transcendental meromorphic solution of the equation \eqref{9.5}, 
    then $\deg_\omega(R)\le 5n+2$, and the following statements hold.
    \begin{itemize}
       \item[\textnormal{(i)}] If $\deg_\omega(Q)=0$, then $\deg_\omega(P)\le n$.
       \item[\textnormal{(ii)}] If $\deg_\omega(Q)\ge 1$, then $\deg_\omega(P)\le \deg_\omega(Q)+n$.
    \end{itemize}
     Moreover, if $Q(z,\omega)$ has a meromorphic function root $b$ in $\omega$ with multiplicity $k$, then $k\le n+1$.
\end{theo}

\begin{remark}
    If $n=1$, then Theorem \ref{the9} reduces into \cite[Theorem 1.1]{n2025}.
\end{remark}

The following example shows the existence of solutions of Theorem \ref{the9}, and $\deg_\omega(P)=\deg_\omega(Q)+n$ is sharp.

\begin{example} 
    The function $\omega(z)=e^{2\pi z}-z$ is a solution of the delay Schwarzian difference equation 
    $$\omega(z+1)-\omega(z-1)+\frac{1}{64\pi^6}(S(\omega,z))^2=\frac{P(z,\omega(z))}{Q(z,\omega(z))},$$
    where
    \begin{align*}
        P(z,\omega(z))&=16\pi^4\omega^6+(96\pi^4-32\pi^3)\omega^5z+(192\pi^4-160\pi^3+25\pi^2)\omega^4z^2+(128\pi^4-224\pi^3  \\
        & \text{ }\text{ }+100\pi^2-6\pi)\omega^3z^3 +(-48\pi^4z^4+168\pi^3z^4+78\pi^2z^4-18\pi z^4+2z^4+2-72\pi^2)\omega^2  \\
        & \text{ }\text{ }+(-96\pi^4+128\pi^3-44\pi^2-18\pi+4)\omega z^5+(4-144\pi^2+24\pi)\omega z+(-32\pi^4+64\pi^3  \\
        & \text{ }\text{ }-47\pi^2+18\pi-1)z^6+(2-72\pi^2+24\pi)z^2-3,  \\
        Q(z,\omega(z))&=16\pi^4\omega^4+(64\pi^4z-32\pi^3)\omega^3+(96\pi^4z^2-96\pi^3z+24\pi^2)\omega^2  \\
        & \text{ }\text{ }+(64\pi^4z^3-96\pi^3z^2+48\pi^2z-8\pi)\omega+(16\pi^4z^4-32\pi^3z^3+24\pi^2z^2-8\pi z+1).
    \end{align*}
   Then we can see that $\deg_\omega(Q)=4$ and $\deg_\omega(P)=6= \deg_\omega(Q)+2$.
\end{example}

The necessary conditions for the existence of subnormal transcendental meromorphic solutions to equation \eqref{eq.1.1} is obtained in the following Theorem \ref{the1}, which can be regarded as the product-theoretic analogue of \cite[Theorem 1.1]{n2025}.

\begin{theo}  \label{the1}
    Let $\omega$ be a subnormal transcendental meromorphic solution of the equation \eqref{eq.1.1}, then $\deg_\omega(R)\le 7$, and the following statements hold.
    \begin{itemize}
       \item[\textnormal{(i)}] If $\deg_\omega(Q)=0$, then $\deg_\omega(P)\le 2$.
       \item[\textnormal{(ii)}] If $\deg_\omega(Q)\ge 1$, then $\deg_\omega(P)\le \deg_\omega(Q)+2$.
    \end{itemize}

Moreover, if $Q(z,\omega)$ has a meromorphic function root $b_1$ in $\omega$ with multiplicity $k$, then $k\le 2$.
\end{theo}

\begin{remark} % where the Remark 2.5
    The condition of subnormal in Theorem \ref{the1} is necessary. For example, it is not difficult to deduce that the function $\omega(z)=e^{e^z}$ is a solution of the delay Schwarzian differential equation
    $$\omega(z+1)\omega(z-1)+2e^zS(\omega,z)=\omega(z)^{e+\frac{1}{e}}-e^{3z}-e^z.$$
    Obviously, $\displaystyle\limsup_{r\rightarrow\infty}\frac{\log T(r,\omega)}{r}=1>0$, here, $Q(z,\omega(z))=1$, $P(z,\omega(z))=\omega(z)^{e+\frac{1}{e}}-e^{3z}-e^z$. Then we have $\deg_\omega(P)=e+\frac{1}{e}>2$ instead of $\deg_{\omega}(P)\le 2$.
\end{remark}

The following two examples show the existence of solutions of case \textnormal{(i)} of Theorem \ref{the1}, and $\deg_\omega(P)=2$ is sharp for transcendental solutions.

\begin{example} \label{exm 1} % where the exm1
    The function $\omega(z)=\tan(\pi z)$ is a solution of the delay Schwarzian difference equation
    $$\omega(z+1)\omega(z-1)+zS(\omega,z)=\omega^2(z)+2\pi^2z,$$
    here $Q(z,\omega(z))=1$, $P(z,\omega(z))=\omega^2(z)+2\pi^2z$. Then we have $\deg_\omega(P)=2$.
\end{example}

\begin{example} \label{exm 2} % where the exm2
    The function $\omega(z)=ze^z$ is a solution of the delay Schwarzian difference equation
    $$\omega(z+1)\omega(z-1)+\frac{2(z+1)^2}{z^2+4z+6}S(\omega,z)=\frac{z^2-1}{z^2}\omega^2(z)-1,$$
    here $Q(z,\omega(z))=1$, $P(z,\omega(z))=\frac{z^2-1}{z^2}\omega^2(z)-1$. Then we have $\deg_\omega(P)=2$.
\end{example}

The following Example \ref{exm 3} and  Example \ref{exm 4} show that the case \textnormal{(ii)} of Theorem \ref{the1} can happen.

\begin{example} \label{exm 3} % where the exm3
    The function $\omega(z)=\frac{1}{e^{\pi z}+1}$ is a solution of the delay Schwarzian difference equation 
    $$\omega(z+1)\omega(z-1)+\frac{-2}{\pi^2}S(\omega,z)=\frac{P(z,\omega(z))}{Q(z,\omega(z))},$$
    where
    \begin{align*}
        P(z,\omega(z))&=(37-24(e^{\pi}+e^{-\pi}))\omega^4+(-84+48(e^{\pi}+e^{-\pi}))\omega^3+(78-29(e^{\pi}+e^{-\pi}))\omega^2  \\
        & \text{ }\text{ }+(-34+5(e^{\pi}+e^{-\pi}))\omega+5,  \\
        Q(z,\omega(z))&=(2-e^{\pi}-e^{-\pi})\omega^2+(-2+e^{\pi}+e^{-\pi})\omega+1  \\
        &=[(2-e^{\pi}-e^{-\pi})\omega-1](\omega-1).
    \end{align*}
   Then we can see that $\deg_\omega(Q)=2$ and $\deg_\omega(P)=4= \deg_\omega(Q)+2$.
\end{example}

\begin{example} \label{exm 4} % where the exm4
    The function $\omega(z)=\frac{1}{\sin(\pi z)+i}$ is a solution of the delay Schwarzian difference equation 
    $$\omega(z+1)\omega(z-1)+S(\omega,z)=\frac{P(z,\omega(z))}{Q(z,\omega(z))},$$
    where
    \begin{align*}
        P(z,\omega(z))&=(4-12\pi^2)\omega^4+(4-12\pi^2)i\omega^3+2\pi^2i\omega,  \\
        Q(z,\omega(z))&=-8i\omega^3+12\omega^2+8i\omega-2  \\
        &=-2(2i\omega-1)(\sqrt{2}\omega+1)(\sqrt{2}\omega-1).
    \end{align*}
   Then we can see that $\deg_\omega(Q)=2$ and $\deg_\omega(P)=4< \deg_\omega(Q)+2=5$.
\end{example}

\section{ Painlev\'{e} \uppercase\expandafter{\romannumeral5} type  difference equation with Schwarzian derivative}

The necessary conditions for the existence of subnormal transcendental meromorphic solutions of equation \eqref{eq.1.3} is presented as follows.
\begin{theo}  \label{the3}%theorem 3.1
    Let $\omega$ be a subnormal transcendental meromorphic solution of \eqref{eq.1.3}, then  $\deg_\omega(R)\le 9$ and  the following statements hold.
    \begin{itemize}
       \item[\textnormal{(i)}] If $\deg_\omega(Q)=0$, then $\deg_\omega(P)\le 4$;
       \item[\textnormal{(ii)}] If $\deg_\omega(Q)\ge 1$, then $\deg_\omega(P)\le \deg_\omega(Q)+4$.
    \end{itemize}
    Moreover, if $Q(z,\omega)$ has a meromorphic function root $b_1$  with multiplicity $k$, then $k\le 2$.

\end{theo}

\begin{remark} % Remark 3.2
    The condition of subnormal in Theorem \ref{the3} is necessary. For example, it is not difficult to deduce that the function $\omega(z)=e^{e^z}$ is a solution of the delay Schwarzian difference equation
    $$(\omega(z)\omega(z+1)-1)(\omega(z)\omega(z-1)-1)+S(\omega,z)=\omega(z)^{2+e+\frac{1}{e}}-\omega(z)^{1+e}-\omega(z)^{1+\frac{1}{e}}-\frac{1}{2}e^{2z}+\frac{1}{2}.$$
    Obviously, $\displaystyle\limsup_{r\rightarrow\infty}\frac{\log T(r,\omega)}{r}=1>0$, here, $Q(z,\omega(z))=1$, $P(z,\omega(z))=\omega(z)^{2+e+\frac{1}{e}}-\omega(z)^{1+e}-\omega(z)^{1+\frac{1}{e}}-\frac{1}{2}e^{2z}+\frac{1}{2}$. Then we have $\deg_\omega(P)=2+e+\frac{1}{e}>4$ instead of $\deg_{\omega}(P)\le 4$.
\end{remark}

The following Example \ref{exm 5} shows the existence of solutions of case $(\textnormal{i})$ of Theorem \ref{the3}.

\begin{example} \label{exm 5} % where the exm5
    The function $\omega(z)=\frac{1}{e^{2\pi iz}}$ is a solution of the delay Schwarzian difference equation  
    $$(\omega(z)\omega(z+1)-1)(\omega(z)\omega(z-1)-1)+\frac{1}{2\pi^2}S(\omega,z)=\omega(z)^4-2\omega(z)^2+2.$$
    Here $Q(z,\omega(z))=1$, $P(z,\omega(z))=\omega(z)^4-2\omega(z)^2+2$. Then we have $\deg_\omega(P)=4$.
\end{example}

The following Example \ref{exm 6} shows the existence of solutions of case $(\textnormal{ii})$ of Theorem \ref{the3}.

\begin{example} \label{exm 6} % where the exm6
    The function $\omega(z)=\frac{1}{e^{z}+1}$ is a solution of the delay Schwarzian difference equation  
    $$(\omega(z)\omega(z+1)-1)(\omega(z)\omega(z-1)-1)+2S(\omega,z)=\frac{P(z,\omega(z))}{Q(z,\omega(z))},$$
    where
    \begin{align*}
        P(z,\omega(z))&=e\omega^4+(e-1)^2\omega^3-(e^2+1)\omega^2,  \\
        Q(z,\omega(z))&=-(e-1)^2\omega^2+(e-1)^2\omega+e.
    \end{align*}
   Then we can see that $\deg_\omega(Q)=2$ and $\deg_\omega(P)=4< \deg_\omega(Q)+4$.
\end{example}

It is worth noting that if $Q(z,\omega)$ has a meromorphic function root $0$  with multiplicity $k$, then  the following Theorem \ref{the4} shows $k\le 1$. Since the proof is similar to that of Theorem \ref{the3}, then we omit its proof.

\begin{theo}  \label{the4}
    Let $\omega$ be a subnormal transcendental meromorphic solution of
    \begin{equation}
    \label{eq.1.4}
        (\omega(z)\omega(z+1)-1)(\omega(z)\omega(z-1)-1)+a(z)S(\omega,z)=\frac{P(z,\omega(z))}{(\omega(z))^k \hat{Q}(z,\omega(z))},
    \end{equation}
    where $k$ is a positive integer, $a(z)$ is small with respect to $\omega$, and  $P(z,\omega)$, $\hat{Q}(z,\omega)$ are polynomials in $\omega$ having meromorphic coefficients small with respect to $\omega$ in $z$, $P(z,\omega)$, $\omega(z)$ and $\hat{Q}(z,\omega)$ are pairwise coprime. Then $k\le 1$.
\end{theo}

\section{Auxiliary results}

The Valiron-Mohon'ko identity is a useful tool to estimate the characteristic function of rational functions. Its proof can be found in \cite[Theorem 2.2.5]{Laine}.

\begin{lem} \cite[Theorem 2.2.5]{Laine} \label{lem1}% where lemma 2.1
Let $\omega $ be a meromorphic function. Then for all irreducible rational functions in $\omega $,
$$R(z, \omega )=\frac{P(z, \omega)}{Q(z, \omega)}=\frac{\sum_{i=0}^{p}a_i(z)\omega ^i}{\sum _{j=0}^qb_j(z)\omega^j},$$
with meromorphic coefficients $a_i(z)$, $b_j(z)$ such that $a_i(z)$ and $b_j(z)$ are small with respect to $\omega$. Then Nevanlinna characteristic function of $R(z, \omega(z))$ satisfies
$$T(r, R(z,\omega))=\deg_{\omega}(R)T(r,\omega)+S(r,\omega),$$
where $\deg_{\omega}(R)=\max \left\{p,q\right\}$ is the degree of $R(z,\omega)$.
\end{lem}

The following two lemmas came from \cite{ZK}. The Lemma \ref{lem2} is the difference version of the logarithmic derivative.

\begin{lem} \cite{ZK}  \label{lem2}% where lemma 2.2
Let $\omega $ be a meromorphic function. If $\omega$ is subnormal, then
$$m(r,\frac{\omega(z+c)}{\omega(z)})=S(r,\omega)$$
hold as $r\notin E$ and $r\rightarrow \infty$, where $E$ is a subset of $[1, +\infty)$ with the zero upper density, that is
$$\overline{dens}E=\limsup_{r\rightarrow \infty}\frac{1}{r}\int_{E\cap \left[1,r\right]}{dt}=0.$$
\end{lem}

\begin{lem} \cite[Lemma 2.1]{ZK}  \label{lem3}
    Let $T(r)$ be a non-decreasing positive function in $[1,+\infty)$ and logarithmic convex with $T(r)\rightarrow +\infty$ as $r\rightarrow\infty$. Assume that
    $$\limsup_{r\rightarrow\infty}\frac{\log T(r)}{r}=0.$$
    Set $\phi(r)=\max_{1\le t\le r}\left\{\frac{t}{\log T(t)}\right\}$. Then given a constant $\delta \in (0,\frac{1}{2})$, we have
    $$T(r)\le T(r+\phi ^{\delta}(r))\le \big(1+4\phi ^{\delta-\frac{1}{2}}(r))T(r), r\notin E_\delta ,$$
    where $E_{\delta}$ is a subset of $[1,+\infty)$ with the zero upper density.
\end{lem}

The following lemma plays an important role in the proof of our results.
\begin{lem} \cite{xc} \label{lem4}% where lemma 2.4
Let $\omega$ be a transcendental meromorphic solution of the delay-differential equation
\begin{eqnarray*}
    &&\varphi(z,\omega )    \\
    &=&\sum_{l\in L}{b_l(z)\omega(z)^{l_{0,0}}\omega(z+c_1)^{l_{1,0}}\cdots \omega(z+c_v)^{l_{v,0}}\left[\omega'(z)\right]^{l_{0,1}}\cdots \left[\omega ^{(\mu)}(z+c_v)\right]^{l_{v,\mu}}}    \\
    &=&0,
\end{eqnarray*}
where $c_1$,..., $c_v$ are distinct complex constants, $L$ is a finite index set consisting of elements of the form $l=(l_{0,0},...,l_{v,\mu})$ and the coefficients $b_l$ are meromorphic functions small with respect to $\omega$ for all $l\in L$. Let $a_1, ..., a_m$ be meromorphic functions small with respect to $\omega$ satisfying $\varphi(z,a_i)\not\equiv 0$ for all $i\in \left\{1,...,m \right\}$. If there exist $s>0$ and $\tau\in (0,1)$ such that
$$\sum_{i=1}^m{n(r,\frac{1}{\omega-a_i})}\le m\tau n(r+s,\omega)+O(1),$$
then
$$\limsup _{r\rightarrow\infty}\frac{\log T(r,\omega )}{r}>0.$$
\end{lem}

The next lemma is a delay-differential version of the Clunie lemma.

\begin{lem} \cite{hw}
    \label{lem5}
    Let $\omega$ be a subnormal transcendental meromorphic solution of 
    \begin{align*}
    P(z,\omega)U(z,\omega)=M(z,\omega),
    \end{align*}
    where $P(z,\omega)$ is a difference polynomial contains just one term of maximal total degree in $\omega$ and its shifts, $U(z,\omega)$ and $M(z,\omega)$ are delay-differential polynomials, if all three with meromorphic coefficients $\alpha_\lambda$ such that $m(r,\alpha_\lambda)=S(r,\omega)$, and the total degree of $\deg(M(z,\omega))\le\deg(P(z,\omega))$. Then
    \begin{align*}
    m(r,U(z,\omega))= S(r,\omega).
    \end{align*}
    \end{lem}

The following lemma is the delay differential version of Mohon'ko theorem, which can be obtained in \cite[Lemma 2.2]{LL}.

\begin{lem} \cite{LL}  \label{lem6}
Let $\omega $ be a subnormal non-rational meromorphic solution of
$$\varphi(r,\omega)=0,$$
where $\varphi(z,\omega )$ is a delay-differential polynomial in $\omega(z)$ with coefficients small with respect to $\omega(z)$. If $\varphi(z,a)\ne 0$ for some small meromorphic function $a(z)$ of $\omega(z)$, then
$$m(r,\frac{1}{\omega-a})=S(r,\omega ).$$
\end{lem}

We next recall some basic properties of the Schwarzian derivative, which can be found in \cite{I1991}.
Let $\omega$ be a meromorphic function.

\begin{itemize}
    \item[\textnormal{(i)}] If $z_0$ is a simple pole of of $\omega(z)$, then $S(\omega,z)$ is regular at $z_0$.
    \item[\textnormal{(ii)}] If $z_0$ is a multiple pole of $\omega(z)$ or a zero of $\omega'(z)$, then $z_0$ is a double pole of $S(\omega,z)$.
\end{itemize}

The following three lemmas are used to estimate the degree of $R(z,\omega)$ in \eqref{9.5}, \eqref{eq.1.1} and \eqref{eq.1.3}, respectively. 

\begin{lem}
\label{l1}
Let $\omega$ be a transcendental meromorphic solution of \eqref{9.5} with $\overline{\lim\limits_{r\to\infty}}\frac{\log T(r,\omega )}{r}$$=0$. Then $\deg_\omega(R)\le 5n+2$. Furthermore, if $\deg_\omega(Q)=0$, then  $\deg_\omega(R)\le n$.
\end{lem}

\noindent{\emph{Proof.}}
By Lemma \ref{lem1}, Lemma \ref{lem2} and Lemma \ref{lem3}, we obtain
\begin{eqnarray*}
    \deg_{\omega }(R)T(r,\omega )+S(r,\omega)&=&T(r,R(z,\omega))    \\
    &\le& T(r,\omega(z+1)-\omega(z-1))+T(r,(S(\omega,z))^n)+S(r,\omega)   \\
    &\le& 2T(r,\omega)+2n\overline{N}(r,\frac{1}{\omega'})+nN(r,\omega)+S(r,\omega)    \\
    &\le& (n+2)T(r,\omega)+2nT(r,\omega')+S(r,\omega)    \\
    &\le& (5n+2)T(r,\omega)+S(r,\omega),
\end{eqnarray*}
which implies that $\deg_{\omega}(R)\le 5n+2$. 

Now, suppose $\deg_\omega(Q)=0$, without loss of generality, assume that $R(z,\omega(z))$ is just $P(z,\omega(z))$. Then \eqref{9.5} becomes
\begin{equation}
    \label{9.6}
       \omega(z+1)-\omega(z-1)+a(z)(S(\omega,z))^n=P(z,\omega(z)).
\end{equation}
Suppose that $\deg_\omega(P)=p\ge n+1$, then $P(z,\omega)=a_p(z)\omega^p(z)+a_{p-1}(z)\omega^{p-1}(z)+\cdots+a_0(z)$, where $a_i(z)$, $i=0,1,...,p$, are meromorphic functions small with respect to $\omega$. We rewrite \eqref{9.6} as
$$\omega(z+1)-\omega(z-1)+a(z)(S(\omega,z))^n-a_{p-1}(z)\omega^{p-1}(z)-\cdots-a_0(z)=a_p(z)\omega^p(z).$$
Applying the Lemma \ref{lem5}, we obtain $m(r,\omega)=S(r,\omega)$, which implies that $N(r,\omega)=T(r,\omega)+S(r,\omega)$. Suppose $z_0$ is a pole of $\omega$ with multiplicity $t$, which is not a zero or pole of the coefficient of \eqref{9.6} and its shift. Then by \eqref{9.6}, at least one term of $\omega(z+1)$ and $\omega(z-1)$ has a pole at $z_0$. Without loss of generality, assuming $\omega(z+1)$ has a pole at $z_0$ with multiplicity $pt$. Then shifting \eqref{9.6} gives
\begin{equation}
    \label{9.7}
    \omega(z+2)-\omega(z)+a(z+1)(S(\omega,z+1))^n=P(z+1,\omega(z+1)).
\end{equation}
It follows from \eqref{9.7} that $\omega(z+2)$ has a pole at $z_0$ with multiplicity $p^2t$,  and $\omega(z+3)$ has a pole at $z_0$  with multiplicity $p^3t$. By continuing the iteration and discussing it in this way, we get
$$n(|z_0|+d,\omega)\ge p^dt+O(1)$$
holds for all positive integer $d$. We have
\begin{equation}
\label{2.4}
    n(\eta, \omega)=\frac{1}{\log \frac{r}{\eta}}\int_{\eta}^r{\frac{dt}{t}}n(\eta,\omega)\le \frac{r}{r-\eta}N(r,\omega)+O(1)\le \frac{r}{r-\eta}T(r,\omega)+O(1)
\end{equation}
for $r>\eta$. Therefore, we can let $r=2|z_0|+2d$ and $\eta=\frac{r}{2}$, then we have
\begin{eqnarray}
    \label{44}
    \limsup_{r\rightarrow\infty}{\frac{\log T(r,\omega)}{r}} &\ge& \limsup_{d\rightarrow\infty}{\frac{\log \left[\frac{1}{2}n(|z_0|+d,\omega)\right]}{2|z_0|+2d}}    \nonumber \\
    &\ge& \limsup_{d\rightarrow\infty}{\frac{\log  (p^dt)-\log 2}{2\left|z_0\right|+2d}}   \nonumber \\
    &=&\frac{\log p}{2}>0.
\end{eqnarray}
This contradicts to the assumption, thus $\deg_\omega(P)\le n$.     $\hfill\Box$

\begin{lem}  \label{lem7}
     Let $\omega$ be a transcendental meromorphic solution of \eqref{eq.1.1} with $\overline{\lim\limits_{r\to\infty}}\frac{\log T(r,\omega )}{r}$$=0$. Then $\deg_\omega(R)\le 7$. Furthermore, if $\deg_\omega(Q)=0$, then  $\deg_\omega(R)\le 2$.
\end{lem}

\noindent{\emph{Proof.}} 
By Lemmas \ref{lem1}, \ref{lem2}, and \ref{lem3}, together with an argument analogous to the proof of Lemma \ref{l1}, we obtain $\deg_{\omega}(R) \le 7$.

Now, suppose $\deg_\omega(Q)=0$, without loss of generality, assume that $R(z,\omega(z))$ is just $P(z,\omega(z))$. Then \eqref{eq.1.1} becomes
\begin{equation}
    \label{2.1}
       \omega(z+1)\omega(z-1)+a(z)S(\omega,z)=P(z,\omega(z)).
\end{equation}
We assume $\deg_\omega(P)=p\ge 3$ and aim to derive a contradiction.
Applying the Lemma \ref{lem5}, we obtain $m(r,\omega)=S(r,\omega)$,  then $N(r,\omega)=T(r,\omega)+S(r,\omega)$.
%$$O(T(r,\omega))=T(r,\omega')+S(r,\omega)=N(r,\omega')\le N(r,\omega)+\overline{N}(r,\omega)\le 2N(r,\omega).$$ 
Suppose $z_0$ is a pole of $\omega$ with multiplicity $m$, which is not a zero or pole of the coefficient of \eqref{2.1} and its shift. Then by \eqref{2.1}, it follows that $z_0$ is a pole of $P(z,\omega(z))$ with multiplicity $mp(\ge 3m)$, therefore, at least one term of $\omega(z+1)$ and $\omega(z-1)$ has a pole at $z_0$. Without loss of generality, assuming $\omega(z+1)$ has a pole at $z_0$ with multiplicity at least $\frac{mp}{2}$. Then shifting \eqref{2.1} gives
\begin{equation}
\label{2.2}
    \omega(z+2)\omega(z)+a(z+1)S(\omega,z+1)=P(z+1,\omega(z+1)).
\end{equation}
It follows from \eqref{2.2} that $\omega(z+2)$ has a pole at $z_0$ with multiplicity at least $\frac{(p^2-2)m}{2}(\ge (\frac{p-\frac{1}{2}}{2})^2m)$. Then shifting \eqref{2.2} gives
\begin{eqnarray*}
    \omega(z+3)\omega(z+1)+a(z+2)S(\omega,z+2)=P(z+2,\omega(z+2)).
\end{eqnarray*}
From the above equation, we have $\omega(z+3)$ has a pole at $z_0$ with multiplicity at least $\frac{(p^3-3p)m}{2}(\ge (\frac{p-\frac{1}{2}}{2})^3m)$. By continuing the iteration and discussing it in this way, we get
$$n(|z_0|+d,\omega)\ge (\frac{p-\frac{1}{2}}{2})^dm+O(1)$$
holds for all positive integer $d$. Since \eqref{2.4},  let $r=2|z_0|+2d$ and $\eta=\frac{r}{2}$. An argument analogous to \eqref{44} likewise yields a contradiction, thus $\deg_\omega(P)\le 2$.    $\hfill\Box$

\begin{lem}  \label{lem8}
     Let $\omega$ be a transcendental meromorphic solution of \eqref{eq.1.3} with $\overline{\lim\limits_{r\to\infty}}\frac{\log T(r,\omega )}{r}$$=0$.Then $\deg_\omega(R)\le 9$. Furthermore, if $\deg_\omega(Q)=0$, then  $\deg_\omega(R)\le 4$.
\end{lem}

\noindent{\emph{Proof.}}
By Lemmas \ref{lem1}, \ref{lem2}, and \ref{lem3}, together with an argument analogous to the proof of Lemma \ref{l1}, we obtain $\deg_{\omega}(R)\le 9$, that is, $\max\left\{\deg_\omega(P),\deg_\omega(Q)\right\}\le 9$.

Now, suppose $\deg_\omega(Q)=0$, without loss of generality, assume that $R(z,\omega(z))$ is just $P(z,\omega(z))$. Then \eqref{eq.1.3} becomes
\begin{equation}
    \label{4.1}
       (\omega(z)\omega(z+1)-1)(\omega(z)\omega(z-1)-1)+a(z)S(\omega,z)=P(z,\omega(z)).
\end{equation}
Suppose that $\deg_\omega(P)=p\ge 5$, \eqref{4.1} can be written as
$$\omega'(z)(\omega(z)\omega(z+1)-1)(\omega(z)\omega(z-1)-1)+a(z)[\omega'''(z)-\frac{3}{2}\frac{\omega''(z)}{\omega'(z)}\omega''(z)]=P(z,\omega(z))\omega'(z).$$
Applying the Lemma \ref{lem5}, we obtain $m(r,\omega')=S(r,\omega)$, which implies that $N(r,\omega')=T(r,\omega')+S(r,\omega)$, then 
$$O(T(r,\omega))=T(r,\omega')+S(r,\omega)=N(r,\omega')\le N(r,\omega)+\overline{N}(r,\omega)\le 2N(r,\omega).$$ 
Suppose $z_1$ is a pole of $\omega$ with multiplicity $t$, which is not a zero or pole of the coefficient of \eqref{4.1} and its shift. Then by \eqref{4.1}, at least one term of $\omega(z+1)$ and $\omega(z-1)$ has a pole at $z_1$. Without loss of generality, assuming $\omega(z+1)$ has a pole at $z_1$ with multiplicity at least $v_1=\frac{(p-2)}{2}t$ $(\ge \frac{p-2}{2})$. Then shifting \eqref{4.1} gives
\begin{equation}
    \label{4.2}
    (\omega(z+1)\omega(z+2)-1)(\omega(z+1)\omega(z)-1)+a(z+1)S(\omega,z+1)=P(z+1,\omega(z+1)).
\end{equation}
It follows from \eqref{4.2} that $\omega(z+2)$ has a pole at $z_1$ with multiplicity at least $v_2=(p-2)v_1-t(> (\frac{p-2}{2})^2)$,  and $\omega(z+3)$ has a pole at $z_1$  with multiplicity $(p-2)v_2-v_1$ $(> (\frac{p-2}{2})^3)$. By continuing the iteration and discussing it in this way, we get
$$n(|z_1|+d,\omega)\ge (\frac{p-2}{2})^d+O(1)$$
holds for all positive integer $d$. Since \eqref{2.4}, let $r=2|z_1|+2d$ and $\eta=\frac{r}{2}$. An argument analogous to \eqref{44} likewise yields a contradiction, thus $\deg_\omega(P)\le 4$.    $\hfill\Box$

\section{Proofs of theorems \ref{the9} and \ref{the1}}

\noindent{\emph{Proof of Theorem \ref{the9}.}} Lemma \ref{l1} shows that $\deg_\omega(R)\le 5n+2$ and the conclusion $(\textnormal{i})$ of Theorem \ref{the9}.

 Next,   we consider the case that $\deg_\omega(Q)=q\ge 1$. We claim $\deg_\omega(P)=p\le \deg_\omega(Q)+n$. Otherwise, supposing $\deg_\omega(P)=p> \deg_\omega(Q)+n$, and aim to derive a contradiction. 
Applying the Lemma \ref{lem5} again, we obtain  $N(r,\omega)=T(r,\omega)+S(r,\omega)$. Suppose $z_0$ is a pole of $\omega$ with multiplicity $t$, which is not a zero or pole of the coefficient of \eqref{9.5} and its shift. Then by \eqref{9.5}, at least one term of $\omega(z+1)$ and $\omega(z-1)$ has a pole at $z_0$. Without loss of generality, assuming $\omega(z+1)$ has a pole at $z_0$. For both $t=1$ and $t\ge2$, $z_0$ is a pole of $\omega(z+1)$ with multiplicity $(p-q)t$. Then shifting \eqref{9.5} gives
\begin{equation}
    \label{9.8}
    \omega(z+2)-\omega(z)+a(z+1)(S(\omega,z+1))^n=\frac{P(z+1,\omega(z+1))}{Q(z+1,\omega(z+1))}.
\end{equation}

We claim that $(p-q)^2>2n$ for $t=1$. Otherwise, supposing $2n\ge (p-q)^2$, then we get $2n\ge (n+1)^2$ since $p-q\ge n+1$, which implies that $n^2+1\le 0$, this is a contradiction. Thus, it follows from \eqref{9.8} that $\omega(z+2)$ has a pole at $z_0$ with multiplicity $(p-q)^2$. If $t\ge 2$, it follows from \eqref{9.8} that $\omega(z+2)$ has a pole at $z_0$ with multiplicity $(p-q)^2t$, and $\omega(z+3)$ has a pole at $z_1$  with multiplicity $(p-q)^3t$. By continuing the iteration and discussing it in this way, we get
$$n(|z_0|+d,\omega)\ge (p-q)^dt+O(1)$$
holds for all positive integer $d$. Since \eqref{2.4}, let $r=2|z_0|+2d$ and $\eta=\frac{r}{2}$, an argument analogous to \eqref{44} likewise yields a contradiction, thus $\deg_\omega(P)\le \deg_{\omega}(Q)+n$. This is the conclusion $(\textnormal{ii})$ of Theorem \ref{the9}.

Finally, we consider the case that the polynomial $Q(z,\omega)$ has a meromorphic function root $b(z)$ small with respect to $\omega$. Then \eqref{9.5} becomes
\begin{equation}
        \label{9.9}
        \omega(z+1)-\omega(z-1)+a(z)(S(\omega,z))^n=\frac{P(z,\omega(z))}{(\omega(z)-b(z))^k\hat{Q}(z,\omega(z))},
\end{equation}
where $P(z,\omega)$, $\omega(z)-b(z)$ and $\hat{Q}(z,\omega)$ are pairwise coprime. Assume that $k\ge n+2\ge 3$, aim for a contradiction. 

Notice that $b$ is not solution of \eqref{9.9}, applying Lemma \ref{lem6}, we obtain $m(r,\frac{1}{\omega-b})=S(r,\omega)$. Thus, $N(r,\frac{1}{\omega-b})=T(r,\omega)+S(r,\omega)$. Hence, we can take one $z_0\in \mathbb{C}$ is a zero of $\omega-b$ with multiplicity $t$,  which is not a zero or pole of the coefficient of \eqref{9.9} and its shift, and $P(z_0,\omega(z_0))\ne 0$. For both $t=1$ and $t\ge2$, it follows from \eqref{9.9} that $z_0$ is a pole of $\omega(z+1)$ with multiplicity at least $kt>2$. Then shifting \eqref{9.9} gives
\begin{equation}
        \label{9.10}
        \omega(z+2)-\omega(z)+a(z+1)(S(\omega,z+1))^n=\frac{P(z+1,\omega(z+1))}{(\omega(z+1)-b(z+1))^k\hat{Q}(z+1,\omega(z+1))}.
\end{equation}

\indent{Case 1.} $\deg_{\omega}(P)\le k+\deg_{\omega}(\hat{Q})$. Then by \eqref{9.10}, we get that $\omega(z+2)$ has a pole at $z_0$ with multiplicity $2n$. Then shifting \eqref{9.10} gives
\begin{equation}
        \label{9.11}
        \omega(z+3)-\omega(z+1)+a(z+2)(S(\omega,z+2))^n=\frac{P(z+2,\omega(z+2))}{(\omega(z+2)-b(z+2))^k\hat{Q}(z+2,\omega(z+2))}.
\end{equation}

 If $t=1$,  then $z_0$ is a pole of $\omega(z+1)$ with multiplicity at least $k>2$, and $z_0$ is a pole of  $(S(\omega,z+2))^n$ with multiplicity $2n$. So it is possible that $z_0+3$ is a zero of $\omega(z)-b(z)$ with multiplicity $t_1$.  By continuing the iteration, if $t_1>1$, then $z_0+4$ is a pole of $\omega$, therefore $t_1=1$.   By considering the multiplicities of  zeros of $\omega-b$ and poles of $\omega$ in the set $\{z_0,z_0+1,z_0+2,z_0+3\}$, we find that there are at least $k+2n$ poles of $\omega$ for $2$ zeros of $\omega-b$. By adding up the contribution from all point $z_0$ to the corresponding counting functions, it follows that 
$$n(r,\frac{1}{\omega-b})\le \frac{2}{k+2n}n(r+3,\omega)+O(1).$$ 
By Lemma \ref{lem4}, we get a contradiction. 

If $t>1$, then $kt>2n$. From \eqref{9.11}, we have $z_0+3$ is is pole of $f$ with multiplicity at least $kt$. By continuing the iteration,  it is possible that $z_0+4$ is a zero of $\omega(z)-b(z)$ with multiplicity $t$, otherwise $z_0+5$ is a pole of  $\omega(z)$. By considering the multiplicities of  zeros of $\omega-b$ and poles of $\omega$ in the set $\{z_0,z_0+1,z_0+2,z_0+3,z_0+4\}$, we find that there are at least $2kt+2n$ poles of $\omega$ for $2t$ zeros of $\omega-b$. it follows that 
$$n(r,\frac{1}{\omega-b_1})\le \frac{t}{kt+n}n(r+4,\omega)+O(1).$$ 
By Lemma \ref{lem4}, we get a contradiction.

\indent{Case 2.} $\deg_{\omega}(P)-k-\deg_{\omega}(\hat{Q})=l\ge1$.  

If $t=1$, then it is possible that $2n=kl$  and $z_0+2$ is a zero of $\omega(z)-b(z)$ with multiplicity $t_2$. By continuing the iteration, we get $t_2=1$, otherwise $z_0+3$ is a pole of $\omega$. As the same as Case 1, we get 
$$n(r,\frac{1}{\omega-b})\le \frac{2}{k}n(r+2,\omega)+O(1).$$ 
Since $k\ge 3$, by Lemma \ref{lem4}, we get a contradiction.

If $t>1$, then by \eqref{9.10}, $z_0+2$ is a pole of $\omega$ 
with multiplicity  $ktl$.  When $l=1$,  it is possible that  and $z_0+3$ is a zero of $\omega(z)-b(z)$ with multiplicity $t_3$. By continuing the iteration, we get $t_3=t$, otherwise $z_0+4$ is a pole of $\omega$. As the same as $t=1$, we get 
$$n(r,\frac{1}{\omega-b})\le \frac{1}{k}n(r+3,\omega)+O(1).$$ 
Since $k\ge 3$, by Lemma \ref{lem4}, we get a contradiction.
When $l>1$, then $z_0+3$ is a pole of $\omega$ with multiplicity $ktl^2$ and $z_0+m$ is a pole of $\omega$ with multiplicity $ktl^m$. Using  the same ideas in Lemma \ref{lem7},  we get a contradiction. Thus, $k\le n+1$.

Therefore, the proof is completed.     $\hfill\Box$

\noindent{\emph{Proof of Theorem \ref{the1}.}} Lemma \ref{lem7} shows that $\deg_\omega(R)\le 7$ and the conclusion $(\textnormal{i})$ of Theorem \ref{the1}.

Next, we consider the case that $\deg_\omega(Q)=q\ge 1$. We claim $\deg_\omega(P)\le \deg(Q)+2$. Otherwise, supposing $\deg_\omega(P)=p\ge \deg_\omega(Q)+3$, and aim to derive a contradiction. 
Applying the Lemma \ref{lem5}, we obtain $m(r,\omega)=S(r,\omega)$.  
Suppose $z_0$ is a pole of $\omega$ with multiplicity $m$, which is not a zero or pole of the coefficient of \eqref{eq.1.1} and its shift. Then by \eqref{eq.1.1}, at least one term of $\omega(z+1)$ and $\omega(z-1)$ has a pole at $z_0$. Without loss of generality, assuming $\omega(z+1)$ has a pole at $z_0$ with multiplicity at least $u_1=\frac{(p-q)}{2}m$ $(\ge \frac{p-q}{2})$. Then shifting \eqref{eq.1.1} gives
\begin{equation}
    \label{3.1}
    \omega(z+2)\omega(z)+a(z+1)S(\omega,z+1)=\frac{P(z+1,\omega(z+1))}{Q(z+1,\omega(z+1))}.
\end{equation}
It follows from \eqref{3.1} that $\omega(z+2)$ has a pole at $z_0$ with multiplicity at least $u_2=(p-q)u_1-m(> (\frac{p-q}{2})^2)$,  and $\omega(z+3)$ has a pole at $z_0$  with multiplicity $(p-q)u_2-u_1$ $(> (\frac{p-q}{2})^3)$. By continuing the iteration and discussing it in this way, we get
$$n(|z_0|+d,\omega)\ge (\frac{p-q}{2})^d+O(1)$$
holds for all positive integer $d$. Since \eqref{2.4}, let $r=2|z_0|+2d$ and $\eta=\frac{r}{2}$, an argument analogous to \eqref{44} likewise yields a contradiction, thus $\deg_\omega(P)\le \deg_{\omega}(Q)+2$. This is the conclusion $(\textnormal{ii})$ of Theorem \ref{the1}.

Finally, we consider the case that  $Q(z,\omega)$ has a meromorphic function root $b_1$  with multiplicity $k$. We rewrite \eqref{eq.1.1} into 
    \begin{equation}
        \label{eq.1.2}
        \omega(z+1)\omega(z-1)+a(z)S(\omega,z)=\frac{P(z,\omega(z))}{(\omega(z)-b_1(z))^k\hat{Q}(z,\omega(z))}.
    \end{equation}
Assuming $k\ge 3$ aims for a contradiction. 
Firstly, assume $b_1(z)\not\equiv 0$. Notice that $b_1$ is not solution of \eqref{eq.1.2}, applying Lemma \ref{lem6}, we obtain $m(r,\frac{1}{\omega-b_1})=S(r,\omega)$. Thus, $N(r,\frac{1}{\omega-b_1})=T(r,\omega)+S(r,\omega)$. Hence, we can take one $z_0\in \mathbb{C}$ is a zero of $\omega-b_1$ with multiplicity $m$,  which is not a zero or pole of the coefficient of \eqref{eq.1.2} and its shift, and $P(z_0,\omega(z_0))\ne 0$. Then by \eqref{eq.1.2}, at least one of $z_0+1$ and $z_0-1$ is a pole of $\omega(z)$. Without loss of generality, suppose that $\omega(z+1)$ has a pole at $z_0$ of multiplicity at least $\frac{km}{2}\ge 2$. Then shifting \eqref{eq.1.2} gives
\begin{equation}
    \label{3.2}
    \omega(z+2)\omega(z)+a(z+1)S(\omega,z+1)=\frac{P(z+1,\omega(z+1))}{(\omega(z+1)-b_1(z+1))^k\hat{Q}(z+1,\omega(z+1))}.
\end{equation}

\indent{Case 1.} $\deg_{\omega}(P)\le k+\deg_{\omega}(\hat{Q})$. Then by \eqref{3.2}, we get that $\omega(z+2)$ has a pole at $z_0$ with multiplicity at least $2$. Then shifting \eqref{3.2} gives
\begin{equation}
    \label{3.3}
    \omega(z+3)\omega(z+1)+a(z+2)S(\omega,z+2)=\frac{P(z+2,\omega(z+2))}{(\omega(z+2)-b_1(z+2))^k\hat{Q}(z+2,\omega(z+2))}.
\end{equation}

\indent{Subcase 1.1.} $\frac{km}{2}=2$, that is, $k=4$, $m=1$. By \eqref{3.3}, we get that $\omega(z+3)$ might have a finite value at $z_0$. Now we consider the most extreme case where $z_0$ is a zero of $\omega(z+3)-b_1(z+3)$ with multiplicity $m_1$. Shifting \eqref{3.3} gives
\begin{equation}
    \label{3.4}
    \omega(z+4)\omega(z+2)+a(z+3)S(\omega,z+3)=\frac{P(z+3,\omega(z+3))}{(\omega(z+3)-b_1(z+3))^k\hat{Q}(z+3,\omega(z+3))}.
\end{equation}

Then by \eqref{3.4}, we get that $\omega(z+4)$ has a pole at $z_0$ with multiplicity $km_1-2\ge 1$. If we continue to iterate \eqref{3.4}, we will eventually obtain the same result as before. Considering all of the zeros of $\omega-b_1$, it follows that 
$$n(r,\frac{1}{\omega-b_1})\le \frac{1}{4}n(r+2,\omega)+O(1).$$
We get a contradiction from Lemma \ref{lem4}.

\indent{Subcase 1.2.} $\frac{km}{2}>2$. By \eqref{3.3}, we get that $\omega(z+3)$ has a zero at $z_0$ with multiplicity $\frac{km}{2}-2$.  Considering all of the zeros of $\omega-b_1$, it follows that 
$$n(r,\frac{1}{\omega-b_1})\le \frac{2m}{km+4}n(r+2,\omega)+O(1).$$
We get a contradiction from Lemma \ref{lem4}.

\indent{Case 2.} $\deg_{\omega}(P)\ge k+\deg_{\omega}(\hat{Q})+1$. Set $\deg_{\omega}(P)-k-\deg_{\omega}(\hat{Q})=n\ge1$. 

\indent{Subcase 2.1.} $\frac{kmn}{2}=2$, that is, $k=4$, $m=1$, $n=1$. By \eqref{3.2}, we get that $\omega(z+2)$ might have a finite value at $z_0$. Now we consider the most extreme case where $z_0$ is a zero of $\omega(z+2)-b_1(z+2)$ with multiplicity $m_2$. If $km_2\le \frac{km}{2}$, that is, $m_2\le\frac{m}{2}$, we get that $m_2=0$ since $m=1$. Considering all of the zeros of $\omega-b_1$, it follows that 
$$n(r,\frac{1}{\omega-b_1})\le \frac{1}{2}n(r+1,\omega)+O(1).$$
We get a contradiction from Lemma \ref{lem4}. If $km_2> \frac{km}{2}$, by \eqref{3.3}, we get that $\omega(z+3)$ has a pole at $z_0$ with multiplicity $km_2- \frac{km}{2}$. If we continue to iterate \eqref{3.3}, we will eventually obtain the same result as before. Considering all of the zeros of $\omega-b_1$, it follows that 
$$n(r,\frac{1}{\omega-b_1})\le \frac{2}{k}n(r+1,\omega)+O(1).$$
We get a contradiction from Lemma \ref{lem4}.

\indent{Subcase 2.2.} $\frac{kmn}{2}>2$. By \eqref{3.2}, we get that $\omega(z+2)$ has a pole at $z_0$ with multiplicity $\frac{kmn}{2}$. If $\frac{kmn^2}{2}=\frac{km}{2}$, that is, $n=1$, then by \eqref{3.3}, we get that $\omega(z+3)$ might have a finite value at $z_0$. Now we consider the most extreme case where $z_0$ is a zero of $\omega(z+3)-b_1(z+3)$ with multiplicity $m_3$. If $km_3\le \frac{kmn}{2}$, that is, $m_3\le\frac{m}{2}$, by \eqref{3.4}, we get that $\omega(z+4)$ might have a finite value at $z_0$. Considering all the zeros of $\omega-b_1$, it follows that 
$$n(r,\frac{1}{\omega-b_1})\le \frac{3}{2k}n(r+3,\omega)+O(1).$$
We get a contradiction from Lemma \ref{lem4}. If $km_3>\frac{kmn}{2}$, we get that $\omega(z+4)$ has a pole at $z_0$ with multiplicity $km_3-\frac{kmn}{2}$, if we continue to iterate \eqref{3.4}, we will eventually obtain the same result as before. Then considering all the zeros of $\omega-b_1$, it follows that 
$$n(r,\frac{1}{\omega-b_1})\le \frac{1}{k}n(r+2,\omega)+O(1).$$
We get a contradiction from Lemma \ref{lem4}.
If $\frac{kmn^2}{2}>\frac{km}{2}$, that is, $n\ge 2$, by \eqref{3.3}, we get that $\omega(z+3)$ has a pole at $z_0$ with multiplicity $\frac{kmn^2}{2}-\frac{km}{2}$, and by \eqref{3.4}, we can get $\omega(z+4)$ has a pole at $z_0$ with multiplicity $\frac{kmn}{2}(n^2-2)$. After \eqref{3.4} iterations, we can obtain that $\omega(z+5)$ has a pole at $z_0$ with multiplicity $\frac{km}{2}(n^4-3n^2+1)$ and so on. Considering all  the zeros of $\omega-b_1$, it follows that 
$$n(r,\frac{1}{\omega-b_1})\le \frac{2}{kn+kn^2}n(r+3,\omega)+O(1).$$
We get a contradiction from Lemma \ref{lem4}.

Now, assume $b_1(z)\equiv 0$. Using the same method as above, we also arrive at a contradiction. Here, we omit the details.

To sum up, $k\le 2$. Therefore, the proof is completed.   $\hfill\Box$

\section{Proof of Theorem \ref{the3}}
Lemma \ref{lem8} shows that $\deg_\omega(R)\le 9$ and the conclusion $(\textnormal{i})$ of Theorem \ref{the3}.

Next, we consider the case that $\deg_\omega(Q)=q\ge 1$. We claim $\deg_\omega(P)-\deg_\omega(Q)\le4$. Otherwise, supposing
$\deg_\omega(P)-\deg_\omega(Q)> 4$, and aim to derive a contradiction. 
Applying the Lemma \ref{lem5}, we obtain $m(r,\omega)=S(r,\omega)$. 
Suppose $z_1$ is a pole of $\omega$ with multiplicity $t$, which is not a zero or pole of the coefficient of \eqref{eq.1.3} and its shift. Then by \eqref{eq.1.3}, at least one term of $\omega(z+1)$ and $\omega(z-1)$ has a pole at $z_1$. Without loss of generality, assuming $\omega(z+1)$ has a pole at $z_1$ with multiplicity at least $\beta_1=\frac{(p-q-2)}{2}t$ $(\ge \frac{p-q-2}{2})$. Then shifting \eqref{eq.1.3} gives
\begin{equation}
    \label{4.3}
    (\omega(z+1)\omega(z+2)-1)(\omega(z+1)\omega(z)-1)+a(z+1)S(\omega,z+1)=\frac{P(z+1,\omega(z+1))}{Q(z+1,\omega(z+1))}.
\end{equation}
It follows from \eqref{4.3} that $\omega(z+2)$ has a pole at $z_1$ with multiplicity at least $\beta_2=(p-q-2)\beta_1-t(> (\frac{p-q-2}{2})^2)$,  and $\omega(z+3)$ has a pole at $z_1$  with multiplicity $(p-q-2)\beta_2-\beta_1$ $(> (\frac{p-q-2}{2})^3)$. By continuing the iteration and discussing it in this way, we get
$$n(|z_1|+d,\omega)\ge (\frac{p-q-2}{2})^d+O(1)$$
holds for all positive integer $d$. Since \eqref{2.4}, let $r=2|z_1|+2d$ and $\eta=\frac{r}{2}$, an argument analogous to \eqref{44} likewise yields a contradiction, thus $\deg_\omega(P)\le \deg_{\omega}(Q)+4$. This is the conclusion $(\textnormal{ii})$ of Theorem \ref{the3}.

Finally, we consider the case that  $Q(z,\omega)$ has a meromorphic function root $b_1$ in $\omega$ with multiplicity $k$. We rewrite \eqref{eq.1.3} into 
     \begin{equation}
    \label{eq.1.5}
        (\omega(z)\omega(z+1)-1)(\omega(z)\omega(z-1)-1)+a(z)S(\omega,z)=\frac{P(z,\omega(z))}{(\omega (z)-b_1(z))^k \hat{Q}(z,\omega(z))},
    \end{equation}
Assuming $k\ge 3$ aims for a contradiction.  Notice that $b_1$ is not solution of \eqref{eq.1.5}, applying Lemma \ref{lem6}, we obtain $m(r,\frac{1}{\omega-b_1})=S(r,\omega)$. Thus, $N(r,\frac{1}{\omega-b_1})=T(r,\omega)+S(r,\omega)$. Hence, we can take one $z_1\in \mathbb{C}$ is a zero of $\omega-b_1$ with multiplicity $t$,  which is not a zero or pole of the coefficient of \eqref{eq.1.5} and its shift, and $P(z_1,\omega(z_1))\ne 0$. Then by \eqref{eq.1.5}, at least one of $z_1+1$ and $z_1-1$ is a pole of $\omega(z)$. Without loss of generality, suppose that $\omega(z+1)$ has a pole at $z_1$ of multiplicity $\frac{kt}{2}\ge2$. Then shifting \eqref{eq.1.5} gives
\begin{eqnarray}
    \label{4.6}
    (\omega(z+1)\omega(z+2)-1)(\omega(z+1)\omega(z)-1)+a(z+1)S(\omega,z+1)  \nonumber \\
    =\frac{P(z+1,\omega(z+1))}{(\omega(z+1)-b_1(z+1))^k\hat{Q}(z+1,\omega(z+1))}.
\end{eqnarray}

\indent{Case 1.} $\deg_{\omega}(P)\le k+\deg_{\omega}(\hat{Q})$. Then by \eqref{4.6}, we get that $\omega(z+2)$ has a zero at $z_1$ with multiplicity $kt-2$. Considering all of the zeros of $\omega-b_1$, it follows that 
$$n(r,\frac{1}{\omega-b_1})\le \frac{2}{k}n(r+1,\omega)+O(1).$$
We get a contradiction from Lemma \ref{lem4}.

\indent{Case 2.} $\deg_{\omega}(P)\ge k+\deg_{\omega}(\hat{Q})+1$. Set $\deg_{\omega}(P)-k-\deg_{\omega}(\hat{Q})=n\ge1$. 

\indent{Subcase 2.1.} $\frac{ktn}{2}\le 2$. By \eqref{4.6}, we get that $\omega(z+2)$ has a zero at $z_1$ with multiplicity $kt-2$. Considering all of the zeros of $\omega-b_1$, it follows that 
$$n(r,\frac{1}{\omega-b_1})\le \frac{2}{k}n(r+1,\omega)+O(1).$$
We get a contradiction from Lemma \ref{lem4}.

\indent{Subcase 2.2.} $\frac{ktn}{2}> 2$. If $n=2$, then by \eqref{4.6}, we get that $\omega(z+2)$ might have a finite value at $z_1$. Now we consider the most extreme case where $z_1$ is a zero of $\omega(z+2)-b_1(z+2)$ with multiplicity $t_2$. Shifting \eqref{4.6} gives
\begin{eqnarray}
    \label{4.7}
    (\omega(z+2)\omega(z+3)-1)(\omega(z+2)\omega(z+1)-1)+a(z+2)S(\omega,z+2)  \nonumber \\
    =\frac{P(z+2,\omega(z+2))}{(\omega(z+2)-b_1(z+2))^k\hat{Q}(z+2,\omega(z+2))}.
\end{eqnarray}

By \eqref{4.7}, we get that $\omega(z+3)$ may have a pole at $z_1$ with multiplicity $kt_2-\frac{kt}{2}$. If we continue to iterate \eqref{4.7}, we will eventually obtain the same result as before. Considering all of the zeros of $\omega-b_1$, it follows that 
$$n(r,\frac{1}{\omega-b_1})\le \frac{2}{k}n(r+1,\omega)+O(1).$$
We get a contradiction from Lemma \ref{lem4}. If $n>2$, it follows from \eqref{4.6} that $\omega(z+2)$ has a pole at $z_1$ with multiplicity $\frac{ktn}{2}-kt$. 

If $2(\frac{ktn}{2}-kt)+\frac{kt}{2}=n(\frac{ktn}{2}-kt)$, that is, $n=3$, then by \eqref{4.7}, we get that $\omega(z+3)$ might have a finite value at $z_1$. Now we consider the most extreme case where $z_1$ is a zero of $\omega(z+3)-b_1(z+3)$ with multiplicity $t_3$. Shifting \eqref{4.7} gives
\begin{eqnarray}
    \label{4.8}
    (\omega(z+3)\omega(z+4)-1)(\omega(z+3)\omega(z+2)-1)+a(z+3)S(\omega,z+3)  \nonumber \\
    =\frac{P(z+3,\omega(z+3))}{(\omega(z+3)-b_1(z+3))^k\hat{Q}(z+3,\omega(z+3))}.
\end{eqnarray}

If $kt_3\le \frac{ktn}{2}-kt$, that is, $t_3\le\frac{t}{2}$, by \eqref{4.8}, we get that $\omega(z+4)$ might have a finite value at $z_1$. Considering all of the zeros of $\omega-b_1$, it follows that 
$$n(r,\frac{1}{\omega-b_1})\le \frac{3}{2k}n(r+3,\omega)+O(1).$$
We get a contradiction from Lemma \ref{lem4}. If $kt_3>\frac{ktn}{2}-kt$, we get that $\omega(z+4)$ has a pole at $z_1$ with multiplicity $kt_3-\frac{ktn}{2}+kt$, if we continue to iterate \eqref{4.8}, we will eventually obtain the same result as before. Then considering all of the zeros of $\omega-b_1$, it follows that 
$$n(r,\frac{1}{\omega-b_1})\le \frac{1}{k}n(r+2,\omega)+O(1).$$
We get a contradiction from Lemma \ref{lem4}.

If $2(\frac{ktn}{2}-kt)+\frac{kt}{2}>n(\frac{ktn}{2}-kt)$. By \eqref{4.7}, we get that $\omega(z+3)$ has a zero at $z_1$ with multiplicity $2(\frac{ktn}{2}-kt)+\frac{kt}{2}-n(\frac{ktn}{2}-kt)$. Considering all of the zeros of $\omega-b_1$, it follows that 
$$n(r,\frac{1}{\omega-b_1})\le \frac{2}{k(n-1)}n(r+2,\omega)+O(1).$$
We get a contradiction from Lemma \ref{lem4}.

If $2(\frac{ktn}{2}-kt)+\frac{kt}{2}<n(\frac{ktn}{2}-kt)$. By \eqref{4.7}, we get that $\omega(z+3)$ has a pole at $z_1$ with multiplicity $n(\frac{ktn}{2}-kt)-2(\frac{ktn}{2}-kt)-\frac{kt}{2}$. Considering all of the zeros of $\omega-b_1$, it follows that 
$$n(r,\frac{1}{\omega-b_1})\le \frac{2}{k(n-1)}n(r+2,\omega)+O(1).$$
We get a contradiction from Lemma \ref{lem4}.

To sum up, $k\le 2$. Therefore, the proof is completed.

\noindent\textbf{Declaration of competing interest}

The authors have no competing interests to declare that are relevant to the content of this article.

%\noindent{\textbf{Discussion:}}\begin{itemize}\item A natural question arises: if we combine the Schwarzian differential equation \eqref{00} with either equation \eqref{1.1} or equation \eqref{09} to form the following equations, $$\omega(z+1)\omega(z-1)+a(z)(S(\omega,z))^n=R(z,\omega(z))$$ or $$(\omega(z)\omega(z+1)-1)(\omega(z)\omega(z-1)-1)+a(z)(S(\omega,z))^n=R(z,\omega(z)),$$would analogous Malmquist-type theorems hold for these two equations?\end{itemize}

%%~~~~~~~~~~~~~~~~~~~~~~~~~~~~~~~~~~~~~~~~~~~~~~~~~reference~~~~~~~~~~~~~~~~~~~~~~~~~~~~~~~~~~~~~~~~~~~~~~~~~~~~~~~~~~~~~~~~~~~~~~~
%%~~~~~~~~~~~~~~~~~~~~~~~~~~~~~~~~~~~~~~~~~~~~~~~~~~~~~~~~~~~~~~~~~~~~~~~~~~~~~~~~~~~~~~~~~~~~~~~~~~~~~~~~~~~~~~

\end{document}